\documentstyle[12pt]{article}
\begin{document}

\title{ Resolution of the Cauchy problem for the Toda lattice \\
with non-stabilized initial data}
\author{ Mikhail Kudryavtsev}
\maketitle

\begin{center}
{\it The Institute for Low Temperature Physics and Engineering\\
of the National Academy of Science of Ukraine,\\
\smallskip
Ukraine, Kharkov, 61103, Lenin Ave., 47}\\
\medskip
E-mail: kudryavtsev@ilt.kharkov.ua
\end{center}
\medskip

\abstract{

This paper is the continuation of the work "On an inverse problem for
finite-difference operators of second order" ([1]). We consider the Cauchy
problem for the Toda lattice in the case when the corresponding $L$-operator
is a Jacobi matrix with bounded elements, whose spectrum of multiplicity $2$
is separated from its simple spectrum and contains an interval of absolutely
continuous spectrum. Using the integral equation of the inverse problem
for this matrix, obtained in the previous work, we solve the Cauchy problem
for the Toda lattice with non-stabilized initial data.

}

\bigskip
\bigskip

\centerline{\bf Introduction}

\bigskip

   {\bf 1.}\  In paper [1] we explained the importance of the extension
of classes of initial data for which the Cauchy problem for nonlinear
evolutionary equations can be solved. Due to this reason the goal of
many investigations is the search of new inverse problems for linear
L-operators, which can be applied to solve the corresponding Cauchy
problems with new and possibly wider classes of initial data.

   In this connection in paper [1] we considered the Cauchy problem
for the equation of oscillation of the doubly-infinite Toda lattice
$$ { d\,^2 \, x_k \over d \, t^2 } =
   e^{x_{k+1} - x_k} - e^{x_k - x_{k-1}} \,, \quad  k\in {\bf Z} , \eqno(0.1)
$$
$$  x_k(0)=v_k, \qquad  \dot x_k (0) = w_k\,. \eqno(0.2)
$$

   The Jacobi matrix
$$
J=\pmatrix{
   \ddots&\ddots&\ddots&{}&{}&{}&{}&{}\cr
   {}&b_{-2}&a_{-1}&b_{-1}&{}&0&{}&{}\cr
   {}&{}&b_{-1}&a_0&b_0&{}&{}&{}\cr
   {}&0&{}&b_0&a_1&b_1&{}&{}\cr
   {}&{}&{}&{}&b_1&a_2&b_2&{}\cr
   {}&{}&{}&{}&{}&\ddots&\ddots&\ddots\cr},
   \eqno(0.3)
$$
with the coefficients $a_k$, $b_k$, defined by the formula
$$ a_k= w_k, \qquad b_{k-1} = {\rm exp} {v_k - v_{k-1} \over 2},  \eqno(0.4)
$$
is the $L$-operator that corresponds to this problem.

   The main result of paper [1] is the finding and solving of a new inverse
spectral problem for Jacobi matrix $J$ when it satisfies the following
main conditions: \
{\it  its elements are bounded, its spectrum of multiplicity $2$ is separated
from its simple spectrum and contains an interval $[a,b]$ of absolutely
continuous spectrum of multiplicity $2$. }
\ (Besides, we imposed some additional technical conditions on the behavior
of the arguments of the Weyl functions of the matrix $J$ in the neighborhood
of its spectrum.)

\medskip

   In the present work we use the obtained inverse problem
(more precisely the integral equation of the inverse problem) to solve the
Cauchy problem (0.1), (0.2) for the Toda lattice in the case when the Jacobi
matrix (0.3), defined from the initial conditions by (0.4),
satisfies the described conditions on the type of spectrum.

\bigskip

   {\bf 2.}\  Before starting the solving of the Cauchy problem let us outline
the main results of paper [1]. We denote by $P_k(\lambda)$, $Q_k(\lambda)$
the solutions of the finite-difference equation
$$b_{k-1}\omega_{k-1} + (a_k - \lambda) \omega_{k} + b_k \omega_{k+1} = 0,
\qquad k\in {\bf Z},  \eqno(0.5)
$$
with initial data
$ \  P_0 (\lambda) = 1, \  P_{-1} (\lambda) = 0, \
    Q_0 (\lambda) =0, \  Q_{-1} (\lambda) =1.
$

   As it is known, for nonreal $\lambda$ equation (0.5) has the Weyl solutions
$$ \varphi^R(k,\lambda) = m^R(\lambda)  P_k(\lambda)
                       - {Q_k(\lambda) \over b_{-1}}, \quad k \in {\bf Z},
$$
$$ \varphi^L(k,\lambda) = - {P_k(\lambda) \over b_{-1}}
                       + m^L(\lambda) Q_k(\lambda), \quad k \in {\bf Z},
$$
such that $\sum_{k=N}^\infty  |\varphi^R(k,\lambda)|^2<\infty$, \
$\sum_{-\infty}^{k=N} |\varphi^L(k,\lambda)|^2<\infty$ for any finite $N$.
Here $m^R(\lambda)$ and $m^L(\lambda)$ are the Weyl functions of the matrix
$J$.

   The functions $m^R(\lambda)$, $m^L(\lambda)$ and the number $b_{-1}$
play the role of spectral data from which the Jacobi matrix $J$ is
reconstructed. Without restriction of generality we
consider the case when $[a,b]=[-2,2]$ and $b_{-1}>0$.\footnote{Namely,
if $[a,b]\ne [-2,2]$, we can consider, instead of $J$ and its Weyl functions
$m^R(\lambda)$, $m^L(\lambda)$, the matrix
$$\tilde J = {4\over b-a} (J - {a+b \over 2} I),
$$
where $I$ is the identity matrix, and its Weyl functions
$$ \tilde m^R(\lambda) = m^R( {a+b \over 2} + {b-a \over 4} \lambda ),
   \qquad
   \tilde m^L(\lambda) = m^L( {a+b \over 2} + {b-a \over 4} \lambda ).
$$
Such $\tilde J$ has the necessary form. So, we can solve the inverse problem,
described in [1], for this new matrix and apply the inverse problem to
solve the Cauchy problem for the Toda lattice with the initial data that
correspond (by (0.4)) to the matrix $\tilde J$. If $\tilde x_n(t)$ is
the solution of the equation of the Toda lattice with these initial data,
then the functions
$$ x_n(t) = \tilde x_n( {b-a \over 4} t ) + 2n\ln({b-a \over 4}) +
   {a+b\over2}t,
$$
as it can be verified, are the solution of the Toda equation with
the original initial data (i.e.~corresponding to the matrix $J$).}

\bigskip

   Instead of the two Weyl functions and two Weyl solutions, which are defined
in the plane of spectral parameter $\lambda$, we introduce one function and
one parameter defined in the $z$-plane, where the variable $z$ and the
parameter $\lambda$ are connected by the relation $z+z^{-1}=\lambda$:
$$n(z)=  \cases{ - b_{-1} m^R(z+z^{-1}), \quad |z|<1, \cr
                 - \frac{1}{ b_{-1} m^L(z+z^{-1})}, \quad |z|>1, \cr }
$$
$$ \psi(k,z) = n(z) P_k(z+z^{-1}) + Q_k(z+z^{-1}) =
\cases{ - b_{-1} \varphi^R(k,z+z^{-1}), \quad |z|<1, \cr
     \frac{ \varphi^L(k,z+z^{-1}) }{ m^L(z+z^{-1}) }, \quad |z|>1. \cr}
   \eqno(0.6)
$$

   The key step of the inverse problem is the suitable factorization of
the function
$${(z-z^{-1})} {\Bigl( n(z)-n(z^{-1}) \Bigr) }^{-1}$$
and the
choice of the factorizing function $R(z)$, which are obtained by the following
theorem:

   {\it The function ${(z-z^{-1})} {\Bigl( n(z)-n(z^{-1}) \Bigr) }^{-1}$
in its domain of holomorphy can be represented in the form of the product of
two functions $R(z)$, $R(z^{-1})$:
$$ { {z-z^{-1}} \over n(z)-n(z^{-1}) } \ = \  R(z) R(z^{-1}).
$$
The function $R(z)$ may only have singularities at such points $z$,
that $z+z^{-1}$ belongs to the spectrum of the matrix $J$. If $z+z^{-1}$
belongs to the absolutely continuous spectrum of multiplicity $2$,
then $R(z)$ may have singularities at both points $z,z^{-1}$, and if
$z+z^{-1}$ belongs to the simple spectrum, then $R(z)$ may only have
singularities at one of the points $z,z^{-1}$.}
(Paper [1] represents the explicit expression for the function $R(z)$.
This theorem is strictly proved in paper [2] in a more general form.)

\medskip

   Further for all $k\in {\bf Z}$ we introduce the functions
$$ g(k,z) = {R(z) \over R(\infty)} z^{-(k+1)} h_k \psi(k,z)\,, \eqno(0.7)
$$
where
$$ h_k = \cases{ b_{-1} \ldots b_{k-1}, \quad   k \ge 0, \cr
                 1,                     \quad   k = - 1, \cr
                 {1\over b_k \ldots b_{-2}},     \quad   k \leq -2. \cr }
   \eqno(0.8)
$$

   The following theorem is the main result of paper [1]:

\medskip

{\bf Theorem.} { \it The function $g(k,z)$ is representable in the form
$$ g(k,z) = 1
   + \int\limits_{-\infty}^\infty
     { u(k,\alpha) s(\alpha) \over 1 - z \alpha^{-1} }
     d \sigma (\alpha)
   + {1\over\pi} \int\limits_{-\pi}^\pi
     { \hat r(e^{i\theta}) u(k,e^{i\theta})
       \over 1 - z e^{-i\theta} } d \theta,   \eqno(0.9)
$$
where the function $u(k,\beta)$ is the solution of the integral equation
$$ \beta^{2(k+1)} u(k,\beta) -
   { \beta \over |\beta| } {m(\beta) \over 2 q(\beta^{-1})} u(k,\beta^{-1}) +
$$
$$ + \, {\rm v.p.}\int\limits_{-\infty}^\infty
     { u(k,\alpha) \over 1 - \beta^{-1} \alpha^{-1}  }
     s(\alpha) d \sigma (\alpha)
   + {1\over\pi} \, {\rm v.p.} \int\limits_{-\pi}^\pi
     { \hat r(e^{i\theta}) u(k,e^{i\theta})
       \over 1 - \beta^{-1} e^{-i\theta}  } d \theta
    = -1,   \eqno(0.10)
$$
and $ s(\alpha) = \cases{1, \quad |\alpha|>1, \cr
                     -1, \quad |\alpha|<1. \cr} $ \
The functions $m(\alpha)$, $q(\alpha)$, the measure $d\sigma(\alpha)$,
defined on the real line, and the function $\hat r(e^{i\theta})$, defined
on the unit circle, are explicitly expressed in terms of the function
$n(z)$ {\rm (see~[1])}.
If the matrix $J$ satisfies the imposed conditions, the equation
(0.10) is uniquely solvable for any $k$ in the class
$L^2 ( {\bf R}_{\sigma} \cup {\bf T}_{\hat m\over\pi} )$.}\,
(Here ${\hat m\over\pi}$ is the Lebesque measure on the unit circle
${\bf T}$, divided by $\pi$.)

   In equation (0.10) the parameter $\beta$ belongs to the unit circle
and the real line. We observe that the support of the measure
$d\sigma(\alpha)$ and the unit circle ${\bf T}$ coincide with
the set of such points $z$, that $z+z^{-1}$ belongs to the spectrum of the
operator $J$.

\medskip

  The collection
$\{ \hat r(e^{i\theta}),{m(\alpha)\over q(\alpha^{-1})},d\sigma(\alpha) \}$,
consisting of the functions
$\hat r(e^{i\theta}),{m(\alpha)\over q(\alpha^{-1})}$ and the measure
$d\sigma(\alpha) $, determines integral equation (3.39) and representation
(3.40). This collection is called the {\it reduced spectral data} of
the Jacobi matrix $J$. Evidently, they play the same role in our inverse
problem as the scattering data in the inverse scattering problem. According
to the last theorem for each $k\in{\bf Z}$, the equations (0.10),
reconstructed according to these data, have a unique solution
in the space $L^2({\bf R}_{\sigma} \cup {\bf T}_{\hat m\over \pi})$.

   In order to reconstruct the Jacobi matrix $J$ according to the given
spectral data it is necessary to solve equation (0.10), then find
the functions $g(k,z)$ using formula (0.9) and then find the elements
$a_k$, $b_k$ of the Jacobi matrix using formulas
$$ b_{k-1}^2 = { g(k,0) \over g(k-1,0) },  \qquad
   a_k = \lim_{z\to\infty} { z(g(k,z)-g(k+1,z)) \over g(k,z)}.  \eqno(0.11)
$$
(The elements $b_k$ can be found up to their sign.) We will prove the last
two formulas later on (see lemma 3).

   The map
$$ J \mapsto n(z) \mapsto \{ \hat r(e^{i\theta}),
             {m(\alpha)\over q(\alpha^{-1})},d\sigma(\alpha) \} \,,
$$
described in the previous sections, is the solution of the direct spectral
problem, and the map
$$\{ \hat r(e^{i\theta}),{m(\alpha)\over q(\alpha^{-1})},d\sigma(\alpha) \}
   \mapsto J
$$
solves the inverse spectral problem. Such are the main results of paper [1].

\bigskip

   {\bf 3.}\  In order to solve the Cauchy problem for the equation of the
Toda lattice we introduce the operator
$$ \hat \Gamma(k,t)= \beta^{2(k+1)} e^{\beta t} I
    + e^{\beta^{-1} t} C_2 \,,    \eqno(0.12)
$$
where $C_2$ is the operator in the space
$L^2 ( {\bf R}_{\sigma} \cup {\bf T}_{\hat m\over\pi} )$,
defined by the formula
$$ (C_2 u) (\beta) = { \beta \over |\beta| }
   {m(\beta) \over 2 q(\beta^{-1}) } u(\beta^{-1})
   + {\rm v.p.}\int\limits_{-\infty}^\infty
     { u(\alpha) \over 1 - \beta^{-1} \alpha^{-1}  }
     s(\alpha) d \sigma (\alpha)
  + {1\over\pi} \, {\rm v.p.} \int\limits_{-\pi}^\pi
     { \hat r(e^{i\theta}) u (e^{i\theta})
       \over 1 - \beta^{-1} e^{-i\theta}  } d \theta \,.    \eqno(0.13)
$$
According to theorem 4 of paper [1], the operator $\hat \Gamma(k,t)$ is always
invertible and the inverse is bounded. In these notations the equation (0.10)
of the inverse problem can be rewritten in the form
$$\hat \Gamma(k,o) u(k,\beta) + 1 = 0 \,.
$$

Let us define in the space
$L^2({\bf R}_{\sigma} \cup {\bf T}_{\hat m\over\pi})$
the operator $P$, which project the whole space into its subspace consisting
of the constant functions:
$$ P w = \hat k \Bigl\{ \int\limits_{-\infty}^\infty \alpha w(\alpha ) s(\alpha ) d \sigma(\alpha )
            + {1\over\pi} \int\limits_{-\pi}^{\pi} e^{i\theta} \hat r(e^{i\theta})
              w(e^{i\theta}) d \theta \Bigr\} ,
      \quad w \in L^2({\bf R}_{\sigma} \cup {\bf T}_{\hat m\over\pi}) \,.
$$
(The coefficient $\hat k \ne 0$ is chosen so that the condition
$ P \, (1) = 1 $ is satisfied.)

   Let us also denote by $N$ the operator of multiplication by
the variable $\beta$ in the space
$L^2({\bf R}_{\sigma} \cup {\bf T}_{\hat m\over\pi})$:
$$ (Nw)(\beta) = \beta w(\beta). $$

\medskip

   The furthest part of the work is organized in the following way.
In section 1 we prove that for arbitrary parameters
$\{ \hat r(e^{i\theta}),{m(\alpha)\over q(\alpha^{-1})},d\sigma(\alpha) \}$,
which guarantee the invertibleness of the operator $\hat \Gamma(k,t)$,
and for the operators $\hat \Gamma(k,t)$, $P$ and $N$, defined as above,
the functions
$$ \tilde x_k (t) = \ln P N^{-1} \hat \Gamma(k,t)^{-1} N^{-1}
                     \hat \Gamma(k+1,t) \, (1) - c_1,
$$
where $c_1$ is a constant number, are the solutions of the Toda equation
(0.1) without the initial condition (0.2). In section 2 we prove that if
we define a Jacobi matrix $J$ according to (0.3), (0.4) from the initial
conditions (0.2), find its reduced spectral data
$\{ \hat r(e^{i\theta}),{m(\alpha)\over q(\alpha^{-1})},d\sigma(\alpha) \}$
and substitute them into the expression of the operators $\hat \Gamma(k,t)$
and $P$, then the last formula gives us the solution of the Cauchy problem
(0.1), (0.2).

  We note that there is a short presentation of this work in paper [3].

\bigskip

   {\bf 4.}\  As it was already said, the main theorems of paper [1] are
proved with some additional conditions for the matrix $J$. Since in what
follows we will use the invertibleness of the operator $\hat \Gamma(k,t)$,
defined from the reduced spectral data, we will also need these conditions
for the main theorem of the present work. So, we rewrite them here. Put
the arguments of the Weyl functions:
$$ \eta^R(\tau)= \lim_{\varepsilon \downarrow 0} \arg m^R(\tau+i\varepsilon), \qquad
   \eta^L(\tau)= \lim_{\varepsilon \downarrow 0} \arg {-1\over m^L(\tau+i\varepsilon)}.
$$

   Let $\tilde\Omega^R$ and $\tilde\Omega^L$ be the sets of singularities
of the functions $m^R(\lambda)$ and ${-1\over m^L(\lambda)}$,
resp.\footnote{That means that $\tilde\Omega^R$ and $\tilde\Omega^L$
are the sets of the points where the functions $m^R(\lambda)$ and
${-1\over m^L(\lambda)}$ are not holomorphic. It is easy to see
that $\tilde\Omega^R$ and $\tilde\Omega^L$ can be defined as the supports
of the measures $d\rho_R(\lambda)$ and $d\rho_L(\lambda)$, where
$\rho_R(\lambda)$, $\rho_L(\lambda)$ are such nondecreasing functions that
$$ b_{-1} m^R(\lambda) = \int_{-\infty}^\infty
     {d\rho_R(\tau) \over \tau - \lambda } \,,  \qquad
   -{1\over b_{-1} m^L(\lambda)} = {\lambda\over b_{-1}} + \beta +
     \int_{-\infty}^\infty {d\rho_L(\tau) \over \tau - \lambda } \,,
$$
with $\beta \in {\bf R}$.}.

\bigskip

   Let, further,
$$ \tilde\Omega_2\ \equiv \ \tilde\Omega^R\cap \tilde\Omega^L,  \qquad
   \tilde\Omega_1\ \equiv \ (\tilde\Omega^R\backslash \tilde\Omega^L) \cup
                 (\tilde\Omega^L\backslash \tilde\Omega^R),
$$
and $\tilde\Omega_2^s \subset \tilde\Omega_2$ be the set of the common
poles of the functions $m^R(\lambda)$ and ${-1\over m^L(\lambda)}$, and
let $\tilde\Omega_2^a \equiv \tilde\Omega_2 \backslash \tilde\Omega_2^s$.

\medskip

We assume that:

A) {\it All the three sets $\tilde\Omega_1$, $\tilde\Omega_2^s$,
$\tilde\Omega_2^a$ have positive mutual distances,
$[-2,2] \subset \tilde\Omega_2^a$, and the set
$\tilde\Omega_2^s$ is finite or empty.}

B) {\it For some $\varepsilon>0$ almost everywhere
(with respect to Lebesque measure) on the set
$\tilde\Omega^a_2$ }
$$ 0<\varepsilon<\eta^R(\alpha)<\pi-\varepsilon,  \quad
   0<\varepsilon<\eta^L(\alpha)<\pi-\varepsilon.
$$

C) {\it In some neighborhood of the set $\tilde\Omega^a_2$
the function $\eta^R(\alpha)-\eta^L(\alpha)$ satisfies the H\"older
condition. }

D) {\it The set $\tilde\Omega^a_2 \backslash [-2,2]$ can be covered
with mutually disjoint intervals $\delta_l$ on each of which the following
inequalities are true: }
$$
{\rm ess}\sup_{\alpha\in\delta_l}\eta^R(\alpha)-
{\rm ess}\inf_{\alpha\in\delta_l}\eta^R(\alpha)<\pi\ , \quad
{\rm ess}\sup_{\alpha\in\delta_l}\eta^L(\alpha)-
{\rm ess}\inf_{\alpha\in\delta_l}\eta^L(\alpha)<\pi\ ,
$$

E) {\it For some small $\varepsilon>0$
and $0<\alpha<\varepsilon$
$$\eta^R(2+\alpha)=\eta^L(2+\alpha)=0,   \qquad
  \eta^R(-2-\alpha)=\eta^L(-2-\alpha)=\pi,
$$
and the functions $\eta^R(\tau)$, $\eta^L(\tau)$ satisfy the
H\"older condition on the interval $[-2,2]$.
}

\bigskip
\bigskip

\begin{center}
{\bf 1. The construction of solutions of the Toda equation \\
        from the integral operators of the form $\hat \Gamma(k,t)$ }
\footnote{The considerations of this section are the application
of more common scheme of constructing of solutions of nonlinear equations
presented in [4] to the case of the Toda lattice and to our form of the
operator $\hat \Gamma(k,t)$}
\end{center}

   We denote the associative ring of the bounded operators in the space
$L^2({\bf R}_{\sigma} \cup {\bf T}_{\hat m\over\pi})$ by
$K(L^2({\bf R}_{\sigma} \cup {\bf T}_{\hat m\over\pi})) \equiv {\bf K}$.
We suppose that these operators depend on the real parameter $t$ and on the
integer parameter $k$. Let us introduce in this ring operators $\partial$
and $\partial_\alpha$ of the differentiation by $t$ and by $k$
$$ \partial \, x(k,t) = x_t(k,t),  \quad
   \partial_\alpha  \, x(k,t) = x(k+1,t) - x(k,t), \ \ \quad
   x(k,t) \in {\bf K},
$$
and an automorphism $\alpha$:
$$ \alpha  (x(k,t)) = x(k+1,t).
$$
It is evident that
$$ \partial_\alpha   = \alpha   - I,
$$
where $I\in{\bf K}$ is a unit operator in
$L^2({\bf R}_{\sigma} \cup {\bf T}_{\hat m\over\pi})$ ($Ix=x$), and that
$$ \alpha  (xy) = \alpha  (x) \alpha  (y).
$$
It is also clear that the operators of differentiation $\partial$ and
$\partial_\alpha$ are commutative. So are the operators $\partial$ and
$\alpha$. We will denote be \ $e$ \ the identical operator in
$L^2({\bf R}_{\sigma} \cup {\bf T}_{\hat m\over\pi})$, which is the unit
element of the ring $\bf K$.

\bigskip

Instead of $\hat \Gamma(k,t)$ we will for now consider the operator
$$ \Gamma(k,t) = \beta^{k+2} e^{\beta t} I
          + \beta^{-k} e^{\beta^{-1} t} C_2
          = \beta^{-k} \hat \Gamma(k,t) \,,      \eqno(1.1)
$$
with the operator $C_2$, defined in (0.13). It is simply to verify immediately
that this operator satisfies the differential equations
$$ \partial \, \Gamma = (\partial_\alpha   + I) \Gamma,   \quad
   \partial^2 \, \Gamma + \Gamma = A (\partial_\alpha  + I)\Gamma,   \eqno(1.2)
$$
where $A \in {\bf K}$ is the operator of multiplication by
$(\beta+\beta^{-1})$ in the space
$L^2({\bf R}_{\sigma} \cup {\bf T}_{\hat m\over\pi})$.
The operator $A$ is a constant operator (i.e.~it does not depend on $t$
and $k$). As it is seen from (1.2), the logarithmic derivative
$$ \gamma = \Gamma^{-1} \partial \Gamma = \Gamma^{-1} (\partial_\alpha   +I) \Gamma = \Gamma^{-1} \alpha  (\Gamma)
$$
of the operator $\Gamma$ is invertible:
$$\gamma^{-1} = \alpha  (\Gamma^{-1}) \Gamma.$$

\medskip

{\bf Lemma 1.}  (see [3])  {\it The logarithmic derivative
$ \gamma = \Gamma^{-1} \partial \Gamma $ of the operator $\Gamma$ satisfies
the equation }
$$ \partial (\gamma^{-1} \partial \gamma) = \gamma^{-1} \alpha  (\gamma) - \alpha  ^{-1} (\gamma^{-1})\gamma \,.
   \eqno(1.3)
$$

\medskip

P r o o f. 1) Let us first suppose that $\gamma$ satisfies the equations
$$ \partial \gamma - \gamma \partial_\alpha (\gamma) = 0,  \qquad
   \partial^2 \gamma + 2 \gamma \partial \gamma - (\gamma^2 + \partial \gamma + e) \partial_\alpha (\gamma) = 0. \eqno(1.4)
$$
Then
$$ \gamma^{-1} \partial \gamma = \partial_\alpha (\gamma),
$$
$$ \gamma^{-1} \partial^2 \gamma + 2 \partial \gamma - (\gamma - \gamma^{-1} \partial \gamma + \gamma^{-1}) \partial_\alpha (\gamma) =0.
$$
Substituting $\gamma^{-1} \partial \gamma$ instead of
$\partial_\alpha (\gamma)$ into the second equality, we have
$$ \gamma^{-1} \partial^2 \gamma + 2 \partial \gamma - \gamma \gamma^{-1} \partial \gamma
   + \gamma^{-1} \partial \gamma \, \gamma^{-1} \partial \gamma
   - \gamma^{-1} \partial_\alpha (\gamma) = 0,
$$
or
$$ \gamma^{-1} \partial^2 \gamma -  \gamma^{-1} \partial \gamma \, \gamma^{-1} \partial \gamma = \gamma^{-1} \partial_\alpha (\gamma) - \partial \gamma.
$$
Thus,
$$ \partial (\gamma^{-1} \partial \gamma) = \partial (\partial_\alpha (\gamma)) = \partial_\alpha \partial \gamma = \alpha (\partial \gamma ) - \partial \gamma.  \eqno(1.5)
$$
On the other hand,
$$ \partial (\gamma^{-1} \partial \gamma) \equiv \gamma^{-1} \partial^2 \gamma -
\gamma^{-1} \partial \gamma \, \gamma^{-1} \partial \gamma
$$
$$ = \gamma^{-1} \partial \gamma \, \gamma^{-1} \partial \gamma + \gamma^{-1}
\partial_\alpha (\gamma) - \partial \gamma - \gamma^{-1} \partial \gamma \, \gamma^{-1} \partial \gamma
   = \gamma^{-1} \partial_\alpha (\gamma) - \partial \gamma,
$$
Hence,
$$ \alpha (\partial \gamma) = \gamma^{-1} \partial_\alpha (\gamma) = \gamma^{-1} \alpha (\gamma) - e,
$$
$$ \partial \gamma = \alpha ^{-1} (\gamma^{-1}) \gamma - e .
$$
Substituting these expressions into the right-hand side of (1.5),
we obtain that $\gamma$ satisfies equation (1.3):
$$  \partial (\gamma^{-1} \partial \gamma)
    = \alpha ( \alpha ^{-1} (\gamma^{-1}) \gamma - e) - (\alpha ^{-1} (\gamma^{-1}) \gamma -e)
    = \gamma^{-1} \alpha  (\gamma) - \alpha  ^{-1} (\gamma^{-1})\gamma \,,
$$
which is the abstract form of the Toda lattice equation.

2) Thus, to prove the lemma, we have to show that
equations (1.2) imply equations (1.4). Let us prove the first of them.
The equation
$$ \partial \Gamma = (\partial_\alpha  + I) \Gamma = \alpha (\Gamma)
$$
is equivalent to
$$ \Gamma (\Gamma^{-1} \partial \Gamma) = \Gamma \Gamma^{-1} (\partial_\alpha  + I) \Gamma,
$$
or
$$ \Gamma \gamma = \Gamma (\Gamma^{-1} \alpha (\Gamma)).
$$
Applying the operator $\partial$ to both sides of the equality, we have
$$ \partial \Gamma \gamma + \Gamma \partial \gamma = \partial \Gamma \, \Gamma^{-1} \alpha (\Gamma) + \Gamma \partial (\Gamma^{-1} \alpha (\Gamma)),
$$
or, multiplying the equality by $\Gamma^{-1}$ from the left-hand side,
$$ \gamma^2 + \partial \gamma = \gamma \Gamma^{-1} \alpha (\Gamma) + \partial (\Gamma^{-1} \alpha (\Gamma))
   = \gamma \Gamma^{-1} \alpha (\Gamma) - \Gamma^{-1} \partial \Gamma \, \Gamma^{-1} \alpha (\Gamma) + \Gamma^{-1} \partial (\alpha (\Gamma))
$$
$$ = \Gamma^{-1} \partial (\alpha (\Gamma)) = \Gamma^{-1} \alpha (\partial \Gamma)
   = \Gamma^{-1} \alpha (\Gamma \gamma)
   = \Gamma^{-1} (\alpha (\Gamma)) \alpha (\gamma) = \gamma \alpha (\gamma) .
$$
Since $\alpha (\gamma) - \gamma = \partial_\alpha (\gamma)$, this has as
a consequence that
$$\partial \gamma - \gamma \partial_\alpha (\gamma) =0,$$
which we, actually, wanted to prove.

Let us now prove the second of the equalities (1.4). The equation
$$ \partial^2 \Gamma + \Gamma = A (\partial_\alpha  + I)\Gamma = A \alpha (\Gamma)
$$
is equivalent to
$$ \Gamma (\Gamma^{-1} \partial^2 \Gamma + e ) = A \Gamma (\Gamma^{-1} \alpha (\Gamma)).
$$
Applying the operator $\partial$ to this equality, we have
$$ \partial \Gamma (\Gamma^{-1} \partial^2 \Gamma + e ) + \Gamma \partial (\Gamma^{-1} \partial^2 \Gamma + e ) =
   A \partial \Gamma (\Gamma^{-1} \alpha (\Gamma)) + A \Gamma \partial (\Gamma^{-1} \alpha (\Gamma)),
$$
or, multiplying it by $\Gamma^{-1}$ from the left-hand side,
$$ \gamma (\Gamma^{-1} \partial^2 \Gamma + e ) + \partial (\Gamma^{-1} \partial^2 \Gamma + e ) =
   \Gamma^{-1} A \Gamma \{ \Gamma^{-1} \partial \Gamma \, \Gamma^{-1} \alpha (\Gamma) + \partial (\Gamma^{-1} \alpha (\Gamma)) \}.
$$
We can eliminate $\Gamma^{-1} A \Gamma$ from the last equality, because
$$ \Gamma^{-1} A \Gamma \cdot \Gamma \alpha (\Gamma) = \Gamma^{-1} A \alpha (\Gamma) =
   \Gamma^{-1} (\partial^2 \Gamma + \Gamma) = \Gamma^{-1} \partial^2 \Gamma + e,
$$
from which
$$ \Gamma^{-1} A \Gamma = (\Gamma^{-1} \partial^2 \Gamma + e) (\Gamma^{-1} \alpha (\Gamma))^{-1}
   = (\Gamma^{-1} \partial^2 \Gamma + e) (\alpha (\Gamma))^{-1} \Gamma.
$$
Thus, the analyzed equation is equivalent to the equality
$$ \gamma (\Gamma^{-1} \partial^2 \Gamma + e ) + \partial (\Gamma^{-1} \partial^2 \Gamma)
$$
$$ = (\Gamma^{-1} \partial^2 \Gamma + e) (\alpha (\Gamma))^{-1} \Gamma
   \Bigl\{ \Gamma^{-1} \partial \Gamma \, \Gamma^{-1} \partial (\Gamma)
   - \Gamma^{-1} \partial \Gamma \, \Gamma^{-1} \alpha (\Gamma)+ \Gamma^{-1} \partial (\alpha (\Gamma)) \Bigr\}
$$
$$ = (\Gamma^{-1} \partial^2 \Gamma + e) (\alpha (\Gamma))^{-1} \Gamma \Gamma^{-1} \alpha (\partial\Gamma)
   = \Bigl[ \partial \Gamma = \Gamma \gamma \Bigr]
$$
$$ = (\Gamma^{-1} \partial^2 \Gamma + e) (\alpha (\Gamma))^{-1} \alpha (\Gamma) \alpha (\gamma)
   = (\Gamma^{-1} \partial^2 \Gamma + e) \alpha (\gamma).
$$

Further, since
$$ \begin{array}{c}
\partial^2 \gamma = \partial^2 (\Gamma^{-1} \partial \Gamma)
   = \partial ( - \Gamma^{-1} \partial \Gamma \, \Gamma^{-1} \partial \gamma + \Gamma^{-1} \partial^2 \Gamma)\\
   = \partial ( - \gamma^2 + \Gamma^{-1} \partial^2 \Gamma )
   = - \gamma \partial \gamma - \partial \gamma \, \gamma + \partial (\Gamma^{-1} \partial^2 \Gamma )
   \end{array}
$$
and
$$ \partial ( \Gamma^{-1} \partial^2 \Gamma ) = - \gamma^2 + \Gamma^{-1} \partial^2 \Gamma ,
$$
we have
$$ \partial ( \Gamma^{-1} \partial \Gamma ) = \partial^2 \gamma
   + \gamma \partial \gamma + \partial \gamma \, \gamma ,
$$
$$ \Gamma^{-1} \partial^2 \Gamma = \partial \gamma + \gamma^2 .
$$
Substituting this in our equality we obtain
$$ \gamma (\partial \gamma + \gamma^2 + e) + \partial^2 \gamma
   + \gamma \partial \gamma + \partial \gamma \, \gamma
   = (\partial \gamma + \gamma^2 + e) ( \partial_\alpha  \gamma + \gamma),
$$
or
$$ \partial^2 \gamma + 2 \gamma \partial \gamma - (\gamma^2 + \partial \gamma + e) \partial_\alpha (\gamma) =0 ,
$$
which is the second of equations (1.4). This proves the lemma.
\hfill\rule{0.5em}{0.5em}\medskip

\medskip

The logarithmic derivative $\gamma (k,t)$ is an operator in the space
$L^2({\bf R}_{\sigma} \cup {\bf T}_{\hat m\over\pi})$. As a function of
$k$ and $t$ it is a solution of nonlinear equation (1.3) is the ring $\bf K$.
Now we are going to construct from this solution a solution of the same
nonlinear equation for scalar functions of $k$ and $t$.
Let $P\in{\bf K}$ is a projector in the space
$L^2(R_{\sigma} \cup T_{\hat m\over\pi})$, i.e.
$P=P^2$, and let it not depend on $k$ and $t$.

\bigskip

{\bf Lemma 2.} \ {\it If the equation
$$ \partial (y^{-1} \partial y) = y^{-1} \alpha ^{-1} (y) - \alpha ^{-1} (y^{-1}) y   \eqno(1.5)
$$
has in the ring $\bf K$ a solution $y=x$, which satisfies the condition
$$ P x = P x P, \eqno(1.6)
$$
then the element $P x P$ is a solution of the same equation in the subring
$P {\bf K} P$. }

\medskip

P r o o f. According to the assumption, the element $P=P^2$ is constant.
Hence, the following identities are true:
$$ P \partial(y) \equiv \partial (Py),  \quad
   P \alpha (y) \equiv \alpha (Py),
$$
$$ \partial(y) P \equiv \partial (yP),  \quad
   \alpha (y) P \equiv \alpha (yP),
$$
from which we have that for the element $x$, satisfying (1.6),
$$ P \partial(x) = \partial (Px) = \partial (PxP) = \partial (PxP)P,
$$
$$ P \alpha ^{\pm 1}(x) = \alpha ^{\pm 1} (Px) = \alpha ^{\pm 1} (PxP)
   = \alpha ^{\pm 1} (PxP) P.
$$
Besides, from $PxP = Px$ it follows that
$$ PxP x^{-1} = P,
$$
and, subsequently, in the subring $P {\bf K} P$ the equality holds:
$$ P x^{-1} = (PxP)^{-1} P.
$$
(Here $PxP$ is scalar. Hence, $(PxP)^{-1}$ makes sense.)
Thus, if the element $x$ satisfies the second equation of (1.3), then
$$ 0 = P (\partial(x^{-1} \partial x) - x^{-1} \alpha (x) \alpha ^{-1} (x^{-1}) x)
   = \partial (P x^{-1} \partial x) - (P x P)^{-1} P \alpha (x) + \alpha^{-1} (Px^{-1}) x
$$
$$ = \partial ( (P x P)^{-1} \partial (PxP)) - (P x P)^{-1} \alpha (PxP)
     + \alpha^{-1} ((PxP)^{-1}) PxP .
$$
\hfill\rule{0.5em}{0.5em}

\medskip

According to the lemma, to construct  the solutions of the Toda equation
(0.1) from $\gamma$, we only need to satisfy the condition (1.6).
Now we remark that the equation (1.5) has a multiplicative group of
transformations, that is, if $\gamma$ satisfies this equations, then so does
the element $N^{-1} \gamma$, where $N$ is an arbitrary constant element.
Hence, if the element $\gamma$ satisfies the equation
$$ \gamma = \gamma P + N (e-P),          \eqno (1.7)
$$
with some constant elements $N$ and $P=P^2$, then
$$ \gamma (e-P) = N(e-P),
$$
$$ N^{-1} \gamma (e-P) = (e-P),
$$
$$ P N^{-1} \gamma (e-P) = 0.
$$
Thus,
$$ P N^{-1} \gamma = P N^{-1} \gamma P,
$$
that is, the element $PN^{-1} \gamma = P N^{-1} \gamma P $ satisfies the
conditions of lemma 2 and it is a solution of equation (1.5) in the
subring $P {\bf K} P$.

Relation (1.7), which guarantees the possibility of projecting, is equivalent
to the following equation for $\Gamma$:
$$ \partial \Gamma (e-P) = \Gamma N (e-P).    \eqno (1.8)
$$
So, if we choose the operator $N$ and the projector $P$ in the space
$L^2({\bf R}_{\sigma} \cup {\bf T}_{\hat m\over\pi})$
so that the operator $\Gamma$ satisfies relation (1.8), then the scalar
function $z(k,t) = P N^{-1} \gamma(k,t) P $ will satisfy the equation
$$ (z(k,t)^{-1} z_t(k,t))_t =
   z(k,t)^{-1} z(k+1,t) - z(k-1,t)^{-1} z(k,t).
$$
In this case the functions $x_k(t) = \ln z(k,t)$ will satisfy equation (0.1)
of the oscillation of the Toda lattice, which is the goal of this section.

   Substituting into (1.8) the expression (1.1) for $\Gamma(k,t)$, we reduce
the equation to the form
$$ (\beta \beta^{2(k+1)} e^{\beta t} I
          + \beta^{-1} e^{\beta^{-1} t} C_2) (I-P) =
   (\beta^{2(k+1)} e^{\beta t} I + e^{\beta^{-1} t} C_2)
   N (I-P).
$$
This equation will hold if the relations
$$ \beta I (I-P) = I N (I-P),
$$
$$ \beta^{-1} C_2 (I-P) = C_2 (I-P).  \eqno(1.9)
$$
are satisfied. Let us define $N\in{\bf K}$ as the operator of multiplication
by the variable $\beta$ in the space
$L^2({\bf R}_{\sigma} \cup {\bf T}_{\hat m\over\pi})$:
$$ (Nw)(\beta) = \beta w(\beta). \eqno(1.10) $$
Then the first of these relations is always fulfilled.
Let us define the projector $P\in{\bf K}$ as
$$ P w =\hat k \Bigl\{ \int\limits_{-\infty}^\infty \alpha w(\alpha ) s(\alpha ) d \sigma(\alpha )
            + {1\over\pi} \int\limits_{-\pi}^{\pi} e^{i\theta} \hat r(e^{i\theta})
              w(e^{i\theta}) d \theta \Bigr\} ,
      \quad w \in L^2({\bf R}_{\sigma} \cup {\bf T}_{\hat m\over\pi})\,.  \eqno(1.11)
$$
(The coefficient $k \ne 0$ is chosen so that the requirement $P^2=P$ is
satisfied, i.e.~so that $ P \, (1) = 1 $. Do not mix the real coefficient
$\hat k$ with the integer index $k$.)
We have to verify the relation
$$ \beta^{-1} C_2 - C_2 N = (\beta^{-1} C_2 - C_2 N) P.
$$
For this let us calculate, according to the definition of $C_2$ in (0.13),
$$ (( \beta^{-1} C_2 - C_2 N ) w ) (\beta)
$$
$$ = \beta^{-1} \, {\rm v.p.}\, \int\limits_{\Omega_0}
       { w(\alpha ) \over 1 - \beta^{-1} \alpha ^{-1} } s(\alpha ) d \sigma(\alpha )
   + {\beta^{-1}\over\pi} \, {\rm v.p.} \, \int\limits_{-\pi}^{\pi}
       { \hat r(e^{i\theta}) w(e^{i\theta}) \over
         1 - \beta^{-1} e^{-i\theta}                        } d \theta
$$
$$ - {\rm v.p.}\, \int\limits_{\Omega_0}
       { w(\alpha ) \alpha \over 1 - \beta^{-1} \alpha ^{-1} } s(\alpha ) d \sigma(\alpha )
   - {1\over\pi} \, {\rm v.p.} \, \int\limits_{-\pi}^{\pi}
       { \hat r(e^{i\theta}) w(e^{i\theta}) e^{i\theta}  \over
         1 - \beta^{-1} e^{-i\theta}                        } d \theta
$$
$$ = {\rm v.p.}\, \int\limits_{\Omega_0}
       { \beta^{-1}  - \alpha  \over 1 - \beta^{-1} \alpha ^{-1} } w(\alpha )s(\alpha ) d \sigma(\alpha )
   - {1\over\pi} \, {\rm v.p.} \, \int\limits_{-\pi}^{\pi}
       { \beta^{-1} - e^{i\theta}  \over
         1 - \beta^{-1} e^{-i\theta}      } \hat r(e^{i\theta}) w(e^{i\theta}) d \theta
$$
$$ = - \int\limits_{\Omega_0}
             \alpha \, w(\alpha ) s(\alpha ) d \sigma(\alpha )
     - {1\over\pi} \int\limits_{-\pi}^{\pi}
             e^{i\theta}  \hat r(e^{i\theta}) w(e^{i\theta}) d \theta
   = - {1 \over\hat k} P w,
$$
because
$$  { \beta^{-1}  - \alpha  \over 1 - \beta^{-1} \alpha ^{-1} }
    = { \alpha ^{-1} \beta^{-1}  - 1  \over 1 - \beta^{-1} \alpha ^{-1} } \alpha
    = - \alpha  \,.
$$
Thus,
$$ \beta^{-1} C_2 - C_2 N = - {1 \over\hat k} P .
$$
Therefore, the relation to verify is equivalent to the equality
$$ - {1 \over\hat k} P = - {1 \over\hat k} P P
$$
and is always trivially satisfied.

\bigskip

{\bf Conclusion: } {\it The functions
$$ x_k(t) = \ln \, (P N^{-1} \Gamma(k,t)^{-1} \Gamma_t(k,t) P)
$$
$$   = \ln \, (P N^{-1} \hat \Gamma(k,t)^{-1} \hat \Gamma_t(k,t) \, (1))
     = \ln \, (P N^{-1} \hat \Gamma(k,t)^{-1} N^{-1} \hat \Gamma (k+1,t)
       \, (1))
\footnote{We remark that since $\hat \Gamma(k,t) = N^k \Gamma(k,t)$, we have
$$ \hat \Gamma(k,t)^{-1} \hat \Gamma_t(k,t)
   = \Gamma(k,t)^{-1} N^{-k} N^k \Gamma_t(k,t)
   = \Gamma(k,t)^{-1} \Gamma_t(k,t).
$$
Besides, it is seen that $\hat \Gamma_t(k,t) = N^{-1} \hat \Gamma(k+1,t)$. }
$$
are solutions of the Toda lattice equation (0.1).
Hence, if $x_k(0)$ and $\dot x_k(0)$ satisfy the initial data, then
$x_k(t)$ will be solutions of the corresponding Cauchy problem for
the Toda lattice equation. So, to solve this Cauchy
problem it is sufficient to choose the Jacobi matrix $J$ in such way that
the functions
$ x_k(t) = \ln \, (P N^{-1} \hat \Gamma(k,t)^{-1} N^{-1}
   \hat \Gamma (k+1,t) \, (1))$,
where $\hat \Gamma(k,t)$ are the operators defined from the inverse problem
for the matrix $J$, satisfy the given initial data when $t=0$.}

   We emphasize that in this section we have never used the specific type
of the measure $d\sigma(\alpha)$ and the functions
$m(\alpha)$, $q(\alpha)$ and $\hat r(e^{i\theta})$. The only thing we needed
was that $\hat \Gamma(k,t)$ has the common form (5.2$^\prime$),
(5.2$^{\prime\prime}$) and is invertible. This means that all the speculations
of this section are true for arbitrary measure
$d\sigma(\alpha)$ and functions $m(\alpha)$, $q(\alpha)$,
$\hat r(e^{i\theta})$, provided they guarantee the invertibleness of
$\hat \Gamma(k,t)$. In particular, we can take an arbitrary Jacobi matrix,
satisfying the conditions A)--E), then find its reduced spectral data
$d\sigma(\alpha)$, $m(\alpha)$, $q(\alpha)$, $\hat r(e^{i\theta})$, then
define the operators $\hat \Gamma(k,t)$ and $P$ from them.
The solutions $x_k(t)$, constructed from these operators, will
be solutions of the Toda lattice equation (0.1) anyway.
But in order to satisfy the given initial data, it is necessary to
find appropriately the matrix $J$.

\bigskip
\bigskip
\centerline{\bf 2. Satisfying if the initial data}
\bigskip

   So, if we have the initial data $v_k$, $w_k$ we define the Jacobi matrix
$J$ of the form (0.3) by (0.4). Let $\hat \Gamma(k,t)$ is the operator
of the inverse problem equation for the matrix $J$.
We define from this operator the solution
$ x_k(t) = \ln \,
  (P N^{-1} \hat \Gamma(k,t)^{-1} N^{-1} \hat \Gamma (k+1,t) P)$
of the Toda lattice equation.

   In this section we are going to prove that $x_k$ satisfy the initial data
$$ x_k(0) = v_k +c_1, \quad \dot x_k(0)=w_k, \eqno(2.1)
$$
where the number $c_1 \in {\bf R}$ does not depend on $k$. Then the solution
of the Cauchy problem will be given by
$$ \tilde x_k(t) = x_k(t) - c_1.
$$
It is evident that to verify (2.1) we have to prove

\bigskip

{\bf Lemma 3.} {\it The functions $x_k(t)$ satisfy the relations}
$$ {\rm exp } (x_k(0) - x_{k-1}(0)) = b_{k-1}^2,   \quad
   \dot x_k(0) = a_k.
$$

\bigskip

P r o o f. \ 1) Let us first see how the elements
$b_k$ and $a_k$ of the matrix $J$ are expressed in terms of $g(k,z)$.
We remind that the Weyl solutions $\psi(k,z)$ satisfy the finite-difference
equation
$$ b_{k-1} \psi(k-1,z) + a_k \psi(k,z) + b_k \psi(k,z)
   = (z+z^{-1}) \psi(k,z)       \eqno(2.2)
$$
We also considered the functions
$$ g(k,z) = z^{-(k+1)}  h_k {R(z)\over R(\infty)}  \psi(k,z).
$$
According to theorem 1 of paper [1], they satisfy the asymptotic formula
$$ g(k,z) \to {R(0)\over R(\infty)} h_k^2,  \quad z \to 0,
$$
from which it follows that
$$ { g(k,0) \over g(k-1,0) } = { h_k^2 \over h_{k-1}^2 } = b_{k-1}^2.
$$
Let us now derive a finite-difference equation for
$g(k,z)$. Multiplying the equation (2.2) by
$z^{-(k+1)} h_k {R(z)\over R(\infty)}$, we obtain
$$ b_{k-1} { h_k \over h_{k-1} } \cdot { z^{-(k+1)} \over z^{-(k-1+1)} }
           g(k-1,z)       +       a_k g(k,z)
   + b_k {h_k \over h_{k+1} } \cdot { z^{-(k+1)} \over z^{-(k+1+1)} }
      g(k+1,z) = (z+z^{-1}) g(k,z)
$$
or
$$ b_{k-1}^2 z^{-1} g(k-1,z) + a_k g(k,z) + z g(k+1,z)
   = (z+z^{-1}) g(k,z) ,
$$
which is equivalent to
$$ z^{-1} (b_{k-1}^2 g(k-1,z)  - g(k,z) ) + a_k g(k,z)
   = z (g(k,z) - g(k+1,z)) .
$$
Letting $z \to \infty$, we have
$$ a_k = \lim_{z\to\infty} { z ( g(k,z) - g(k+1,z) ) \over g(k,z) }.
$$

2) Using formula (0.9)
$$ g(k,z) = 1
   + \int\limits_{\Omega_0}
     { u(k,\alpha) s(\alpha) \over 1 - z \alpha^{-1} }
     d \sigma (\alpha)
   + {1\over\pi} {\rm v.p.} \int\limits_{-\pi}^\pi
     { \hat r(e^{i\theta}) u(k,e^{i\theta})
       \over 1 - z e^{-i\theta} } d \theta  \,,
$$
we find that
$$ b_{k-1}^2 =
  { 1 + { \displaystyle \int\limits_{\Omega_0} }
     u(k,\alpha) s(\alpha) d \sigma (\alpha)
   + { \displaystyle {1\over\pi} \int\limits_{-\pi}^\pi }
     \hat r(e^{i\theta}) u(k,e^{i\theta}) d \theta
                                                          \over
     1 + { \displaystyle \int\limits_{\Omega_0} }
     u(k-1,\alpha) s(\alpha) d \sigma (\alpha)
   + {1\over\pi} { \displaystyle \int\limits_{-\pi}^\pi }
     \hat r(e^{i\theta}) u(k-1,e^{i\theta}) d \theta  }
   = {  1 + \hat P N^{-1} u_k   \over
        1 + \hat P N^{-1} u_{k-1}           } ,
$$
where $\hat P= {1\over\hat k} P$ (i.e. it is the "non-normalized" projector),
and $u_k \equiv u(k,\beta)$ is the solution of the inverse problem equation
(for $t=0$). We also see that
$$ a_k = \lim_{z\to\infty} z (g(k,z) - g(k+1,z))
$$
$$ \begin{array}{r}
= \lim_{z\to\infty}
  \Bigl\{ \int\limits_{\Omega_0}
     { z \over 1 - z \alpha^{-1} } s(\alpha) (u(k,\alpha) - u(k+1,\alpha))
     d \sigma (\alpha)\\
   + {1\over\pi} \int\limits_{-\pi}^\pi
     { z \over 1 - z e^{-i\theta} } \hat r(e^{i\theta})
     (u(k,e^{i\theta}) - u(k+1,e^{i\theta})) d \theta  \Bigr\}
   \end{array}
$$
$$ = \Bigl\{  - \int\limits_{\Omega_0}
       \alpha s(\alpha) (u(k,\alpha) - u(k+1,\alpha))
       d \sigma (\alpha)
     + {1\over\pi} \int\limits_{-\pi}^\pi
       e^{i\theta} \hat r(e^{i\theta})
       (u(k,e^{i\theta}) - u(k+1,e^{i\theta})) d \theta  \Bigr\}
$$
$$ = \Bigl\{ - \hat P u_k + \hat P u_{k+1} \Bigr\}.
$$

3) Thus, we have to verify that
$$ {  1 + \hat P N^{-1} u_k   \over
      1 + \hat P N^{-1} u_{k-1}       }
   = e^{ x_k - x_{k-1} }
$$
and
$$ \dot x_k =  \Bigl\{ \hat P u_{k+1} - \hat P u_k \Bigr\}.
$$
Let us check the first of these equalities.
$$ b_{k-1} \equiv {  1 + \hat P N^{-1} u_k   \over
                1 + \hat P N^{-1} u_{k-1}       }
     = {\hat k \hat P 1 + \hat P N^{-1} \hat \Gamma_k^{-1} (-1) \over
        \hat k \hat P 1 + \hat P N^{-1} \hat \Gamma_{k-1}^{-1} (-1)  },
$$
where we denoted, for the sake of brevity, $\hat \Gamma_k = \hat \Gamma(k,0)$.
We remind that $\hat \Gamma_k u_k \equiv \hat \Gamma(k,0) u(k,\beta) = -1$.

In the right-hand side of the obtained equality let us express $\hat k \hat P 1$
in the form 
$$ \hat P N^{-1} \hat \Gamma_k^{-1} N^{-1} \tilde \omega_k,$$ 
where
$\tilde \omega_k \in L^2({\bf R}_{\sigma} \cup {\bf T}_{\hat m\over\pi}).$
For this it is sufficient to solve the following equation for a unknown
function $\tilde \omega_k$:
$$\hat k \, 1 = N^{-1} \hat \Gamma_k^{-1} N^{-1} \tilde \omega_k .
$$
The solution of this equation is the function
$$ \tilde \omega_k = N \hat \Gamma_k N \cdot\hat k 1 =\hat k \beta^{2(k+2)} 1
         +\hat k \beta C_2 N 1.
$$
Besides, we remark now that $-1 = - N^{-1} \cdot \beta 1$. Hence,
$$ b_{k-1} = {  \hat P N^{-1} \hat \Gamma_k^{-1} N^{-1}
                  \{\hat k \beta^{2(k+2)} 1  +
                 \hat k \beta C_2 N 1  -  \beta 1 \} \over
                \hat P N^{-1} \hat \Gamma_{k-1}^{-1} N^{-1}
                \{\hat k \beta^{2(k+1)} 1  +
               \hat k \beta C_2 N 1  -  \beta 1 \}       } \,.
$$
To simplify this expression, we show that
$$ \hat \Gamma_{k+1} (k 1) =\hat k \beta^{2(k+2)} 1
       + \hat k \beta C_2 N 1  -  \beta 1 ,  \qquad
 \hat \Gamma_k (k 1) =\hat k \beta^{2(k+1)} 1
       + \hat k \beta C_2 N 1  -  \beta 1 ,
$$
which is equivalent (taking into account the specific of $\hat \Gamma_k$)
to
$$\hat k \beta C_2 N 1  -  \beta 1  =  C_2 (\hat k 1),
$$
or
$$\hat k (C_2 N - \beta^{-1} C_2) \cdot 1 =1.   \eqno (2.3)
$$
But, according to equality (1.9),
$C_2 N - \beta^{-1} C_2 =  \hat P$. Further, since
$\hat k \hat P 1 = 1$, equality (2.3) always holds and
$$ b_{k-1} = {  \hat P N^{-1} \hat \Gamma_k^{-1} N^{-1}
                  \{ \hat \Gamma_{k+1} (k ) \}  \over
                \hat P N^{-1} \hat \Gamma_{k-1}^{-1} N^{-1}
                  \{ \hat \Gamma_k (k ) \}  }
   = { P N^{-1} \hat \Gamma_k^{-1} N^{-1} \hat \Gamma_{k+1} \cdot 1  \over
       P N^{-1} \hat \Gamma_{k-1}^{-1} N^{-1} \hat \Gamma_k \cdot 1  }
   = { P N^{-1} \hat \Gamma_k^{-1} N^{-1} \hat \Gamma_{k+1} P \over
       P N^{-1} \hat \Gamma_{k-1}^{-1} N^{-1} \hat \Gamma_k P }
$$
$$ = { P N^{-1} \gamma_k P \over
       P N^{-1} \gamma_{k-1} P }
   = e^{x_k-x_{k-1}}
$$
(with $\gamma_k \equiv \gamma(k,0)$), what we needed.
Let us now verify that $\dot x_k = a_k$. First,
$ (\hat \Gamma_k^{-1})_t =
  - \hat \Gamma_k^{-1} N^{-1} \hat \Gamma_{k+1} \hat \Gamma_k^{-1} $.
So,
$$ (\gamma_k)_t = (\hat \Gamma_k^{-1} N^{-1} \hat \Gamma_{k+1})_t
       = - \hat \Gamma_k^{-1} N^{-1} \hat \Gamma_{k+1} \hat \Gamma_k^{-1} N^{-1} \hat \Gamma_{k+1}
         + \hat \Gamma_k^{-1} N^{-2} \hat \Gamma_{k+1}.
$$
Hence,
$$ \dot x_k \equiv (\ln PN^{-1} \gamma_k P)_t
   = { P N^{-1} \gamma_k P \over  P N^{-1} \gamma_k P  }
   = { P N^{-1} \gamma_k \cdot 1   \over   P N^{-1} \gamma_k \cdot 1  }
$$
$$ = { P N^{-1} \{ - \hat \Gamma_k^{-1} N^{-1} \hat \Gamma_{k+1} \hat \Gamma_k^{-1} N^{-1} \hat \Gamma_{k+1}
                   + \hat \Gamma_k^{-1} N^{-2} \hat \Gamma_{k+2}  \} \cdot 1
       \over   P N^{-1} \hat \Gamma_k^{-1} N^{-1} \hat \Gamma_{k+1} \cdot 1  }.
$$
Calculate separately the braces $\{\}$ in the nominator:
$$ (I) \equiv \{ - \hat \Gamma_k^{-1} N^{-1} \hat \Gamma_{k+1} \hat \Gamma_k^{-1} N^{-1} \hat \Gamma_{k+1}
                 + \hat \Gamma_k^{-1} N^{-2} \hat \Gamma_{k+2}  \} \cdot 1 =
$$
$$ = - \hat \Gamma_k^{-1} N^{-1} \hat \Gamma_{k+1} \hat \Gamma_k^{-1}
       \{  \beta^{2(k+2)-1}  +  \beta^{-1} C_2 \cdot 1 \}
     + \hat \Gamma_k^{-1} N^{-1} \{  \beta^{2(k+3)-1}  +  \beta^{-1} C_2 \cdot 1 \}\,.
$$
But $\beta^{-1} C_2 \cdot 1 - C_2 N \cdot 1 = -\hat k^{-1}$, from where
$$ \beta^{2(k+2)-1}  +  \beta^{-1} C_2 \cdot 1
   = \beta^{2(k+2)-1}  +  C_2 \cdot \beta  -\hat k^{-1}
   = \hat \Gamma_k \cdot \beta -\hat k^{-1}.
$$
Thus,
$$ (I) = - \hat \Gamma_k^{-1} N^{-1} \hat \Gamma_{k+1} \hat \Gamma_k^{-1}
       \Bigl\{  \hat \Gamma_k \cdot \beta -\hat k^{-1} \Bigr\}
     + \hat \Gamma_k^{-1} N^{-1} \Bigl\{  \hat \Gamma_{k+1} \cdot \beta -\hat k^{-1} \Bigr\}
$$
$$ = \hat \Gamma_k^{-1} N^{-1} \hat \Gamma_{k+1}
        \Bigl\{ - \beta + \hat \Gamma_k^{-1}\hat k^{-1} \Bigr\}
     + \hat \Gamma_k^{-1} N^{-1} \Bigl\{  \hat \Gamma_{k+1} \cdot \beta -\hat k^{-1} \Bigr\}\,.
$$
But $\hat \Gamma_k u_k = - 1$, so
$ -\hat k^{-1} \hat \Gamma_k u_k  = \hat k^{-1}  $  and
$ \hat \Gamma_k^{-1}\hat k^{-1} =\hat k^{-1} u_k $. Hence,
$$ (I) = \hat \Gamma_k^{-1} N^{-1} \hat \Gamma_{k+1}
       \{  - \beta -\hat k^{-1}  u_k \}
     + \hat \Gamma_k^{-1} N^{-1} \{  \hat \Gamma_{k+1} \beta -\hat k^{-1} \}
$$
$$ = \hat \Gamma_k^{-1} N^{-1}
       \Bigl\{ - \hat \Gamma_{k+1} \beta -\hat k^{-1} \hat \Gamma_{k+1} u_k
          + \hat \Gamma_{k+1} \beta -\hat k^{-1} \Bigr\}
   = - {1\over\hat k} \hat \Gamma_k^{-1} N^{-1}
       \Bigl\{ 1 + \hat \Gamma_{k+1} u_k \Bigr\}.
$$
Let us calculate the denominator in the expression for $x_k$:
$$ P N^{-1} \hat \Gamma_k^{-1} N^{-1} \hat \Gamma_{k+1} \cdot 1
   = P N^{-1} \hat \Gamma_k^{-1} N^{-1} \hat \Gamma_{k+1} \cdot {P u_k \over P u_k}
$$
$$ = {1\over P u_k } P N^{-1} \hat \Gamma_k^{-1} N^{-1} \hat \Gamma_{k+1}  P u_k
   = {1\over P u_k } P N^{-1} \hat \Gamma_k^{-1} N^{-1} \hat \Gamma_{k+1} \cdot u_k .
$$
Thus,
$$ \begin{array}{r}
   \dot x_k \equiv { P N^{-1} \gamma_k P \over  P N^{-1} \gamma_k P  }
   = { - { \displaystyle 1\over\hat k} \hat P N^{-1} \Gamma_k^{-1} N^{-1}
       \{ 1 + \hat \Gamma_{k+1} u_k \}   \over
   { \displaystyle 1\over P u_k } P N^{-1} \hat \Gamma_k^{-1} N^{-1} \hat \Gamma_{k+1} \cdot u_k } \\
     =  - {P u_k \over\hat k}
     \Bigl\{ 1 + { P N^{-1} \hat \Gamma_k^{-1} N^{-1} \cdot 1   \over
                 P N^{-1} \hat \Gamma_k^{-1} N^{-1} \hat \Gamma_{k+1} \cdot u_k }
     \Bigr\}.
   \end{array}
$$
We calculate separately the second term:
$$  { P N^{-1} \hat \Gamma_k^{-1} N^{-1} \cdot 1   \over
      P N^{-1} \hat \Gamma_k^{-1} N^{-1} \hat \Gamma_{k+1} \cdot u_k }
  = { P N^{-1} \hat \Gamma_k^{-1} N^{-1} \cdot 1   \over
      { \displaystyle P u_k  \over  P u_{k+1} }
       P N^{-1} \hat \Gamma_k^{-1} N^{-1} \hat \Gamma_{k+1} P u_k }
$$
$$ = {P u_{k+1} \over P u_k }   \cdot
     { P N^{-1} \hat \Gamma_k^{-1} N^{-1} \cdot 1   \over
       P N^{-1} \hat \Gamma_k^{-1} N^{-1} (-1)  }
   = - {P u_{k+1} \over P u_k }.
$$
Finally, substituting this into the expression for $x_k$, we obtain
$$ x_k =  - {P u_k \over\hat k}
     \Bigl\{ 1 - {P u_{k+1} \over P u_k }
     \Bigr\}
  =  - {1 \over\hat k }
     \Bigl\{ P u_{k+1} - P u_k  \Bigr\}
   =  - {1 \over\hat k }
     \Bigl\{ \hat P u_{k+1} - \hat P u_k  \Bigr\}   =  a_k,
$$
which ends the proof of lemma 3.
\hfill\rule{0.5em}{0.5em}

\medskip

In this section we proved the main

\bigskip

{\bf Theorem \ (the main theorem)} \ (on the existence of the solution
of the Cauchy problem for the Toda lattice). \quad
{\it Let the Jacobi matrix $J$, defined by formulas (0.3), (0.4),
satisfies conditions A)-E).

  Then problem (0.1), (0.2) has a solution
$$ \tilde x_k (t) =
   \ln P N^{-1} \hat \Gamma^{-1}(k,t) N^{-1} \hat \Gamma(k+1,t) \, (1) - c_1 \,,
   \eqno(2.4)
$$
where the operators $P$, $N$, $\hat \Gamma (k,t)$ in the space
$L^2({\bf R}_{\sigma} \cup {\bf T}_{\hat m\over\pi})$,
are defined by formulas (1.11), (1.10), (0.12), and $c_1$ is a constant
number.}
\footnote{In order to make clearer the way of finding $\tilde x_k(t)$
from the operators $\hat \Gamma(k,t)$, $P$ ¨ $N$, we comment the last formula.
Evidently, when the operator $\hat \Gamma(k+1,t)$ operate the constant
function 1, we obtain another function in the space
$L^2({\bf R}_{\sigma} \cup {\bf T}_{\hat m\over\pi})$. (Not only is not it
constant on the variable $\beta$ of the space, but also it depends on
$k,t$.) The operator $N^{-1}$ multiplies it by $\beta^{-1}$. The operator
$\hat \Gamma^{-1}(k,t)$ makes of it a new function in the space
$L^2({\bf R}_{\sigma} \cup {\bf T}_{\hat m\over\pi})$, depending on
$k$ and $t$; $N^{-1}$ multiplies the new function by $\beta^{-1}$.
Finally, the projector $P$ makes of the function a constant in the space
$L^2({\bf R}_{\sigma} \cup {\bf T}_{\hat m\over\pi})$. However, it
depends on $k$ and $t$, as before. Taking the logarithm
of this constant, we find
$ \tilde x_k (t) =
   \ln P N^{-1} \hat \Gamma_k^{-1} N^{-1} \hat \Gamma_{k+1} \, (1) - c_1 $.
}

\medskip

The constant number $c_1$, evidently, can be found from the formula
$$ c_1 = \tilde x_0 (0) -
   \ln P N^{-1} \hat \Gamma^{-1}(0,0) N^{-1} \hat \Gamma(1,0) \, (1).
$$

\bigskip

{\bf Remark.}\  The main result can be exposed in more usual terms of the
evolution of the spectral data. The inverse problem equation with spectral
data depending on time, is rewritten in the form

$$ \beta^{2(k+1)} u(k,\beta) -
   { \beta \over |\beta| } {1 \over \tilde r(\beta,t)} u(k,\beta^{-1}) +
$$
$$ + \, {\rm v.p.}\int\limits_{-\infty}^\infty
     { u(k,\alpha) \over 1 - \beta^{-1} \alpha^{-1}  }
     s(\alpha) d \sigma (\alpha,t)
   + {1\over\pi} \, {\rm v.p.} \int\limits_{-\pi}^\pi
     { \hat r(e^{i\theta},t) u(k,e^{i\theta})
       \over 1 - \beta^{-1} e^{-i\theta}  } d \theta
    = e^{-\beta^{-1}t},   \eqno(2.5)
$$
where
$$   \tilde r(\beta,0) = {2q(\beta^{-1}\over m(\beta)}\,, \quad
     d \sigma (\alpha,0) = d \sigma (\alpha)\,, \qquad
     \hat r(e^{i\theta},0) = \hat r(e^{i\theta}) \eqno(2.6)
$$
are the reduced spectral data of the matrix $J$, defined by the initial data.

It is evident that for $t=0$ equation (2.5) coincide with (0.10).

   The time evolution of the reduced spectral data
is determined by the formula
$$   \tilde r(\beta,t) = e^{(\beta-\beta^{-1})t} \tilde r(\beta,0)
     = {2q(\beta^{-1}\over m(\beta)}\,, \eqno(2.7')
$$
$$   d \sigma (\alpha,t) = e^{(\beta-\beta^{-1})t} d \sigma (\alpha,0)
     = e^{(\beta-\beta^{-1})t} d \sigma (\alpha)\,,  \eqno(2.7'')
$$
$$   \hat r(e^{i\theta},t) = e^{(\beta-\beta^{-1})t} \hat r(e^{i\theta},0)
     = e^{(\beta-\beta^{-1})t} \hat r(e^{i\theta})\,. \eqno(2.7''')
$$

\medskip

{\it In order to find the solution of the Cauchy problem
for the Toda lattice (0.1), (0.2) at the time $t$, one need:

1) to find the reduced spectral data
$\{ {2q(\alpha^{-1})\over m(\beta)},d\sigma(\alpha),\hat r(e^{i\theta}) \}$
of the matrix $J$, defined by (0.3), (0.4);

2) multiplying the reduced spectral data by
$e^{(\beta-\beta^{-1})t}$, to find the data
$\{ \tilde r(\beta,t)$, $d\sigma (\alpha,t)$, $\hat r(e^{i\theta},t) \}$;

3) to solve equation (2.5) for given $t$ and parameters
$\tilde r(\beta,t), d\sigma (\alpha,t), \hat r(e^{i\theta},t)$, with right-hand
side $e^{-\beta^{-1}t}$;

4) to reconstruct according to formulas (0.9), (0.11) the matrix $J(t)$ from
the solutions and to find
$$ \tilde x_k(t) = \tilde x_0(t) + 2 \sum_{j=1}^k \ln b_{j-1}(t) \,.
$$

}

\bigskip
\bigskip
\centerline{\bf Conclusions}
\bigskip

   A new type of inverse problem is introduced for the Jacobi matrices
with bounded elements, whose spectrum of multiplicity $2$ s separated from
the simple spectrum and contains an interval of absolutely continuous
spectrum. The spectral data in this problem are explicitly expressed
in terms of the Weyl functions $m^R(\lambda)$, $m^L(\lambda)$ of the
Jacobi matrix on the right and left semiaxis. These spectral data play the
same role in the inverse problem as the scattering data in the classical
inverse scattering problem. The uniquely solvable integral equation allowing
to reconstruct the matrix by these data is obtained.

   The Jacobi matrices that satisfy the condition on the spectrum written
above, are the $L$-operators of the equation of the oscillation of the Toda
lattice with quite a wide class of initial data, which are not stabilized.
(This class includes the already investigated cases of fast-stabilized and
periodic initial data.) This allowed to apply the inverse problem integral
equation obtained to solve the Cauchy problem for the Toda lattice with
non-stabilized initial data.

\bigskip

Acknowledgments. \ {\it This work is partially supported by
INTAS 2000-272. }

\newpage

\centerline{\Large REFERENCES}

\bigskip

\begin{flushleft}
1. Kudryavtsev M. \, On an inverse problem for finite-difference operators
of second order. http://arXiv.org/abs/math.SP/0110276 \\
2. Kudryavtsev M.A. \, The Riemann problem with additional singularities
(in Russian). Mat. fiz., analiz, geom. (2000) v.~7, ü~2.
p.~196--208. English translation: http://arXiv.org/abs/math.SP/0110242 \\
3. Kudryavtsev M. \,  Solution of the Cauchy problem for a Toda lattice
with initial data that are not stabilized. Dokl. NAN Ukrainy.
(2001), ü~3, p.~14-19. \\
4. Marchenko V.A. \, Nonlinear equations and operator algebras. Kiev:
"Naukova dumka" (1986)\\
5. G.~Teschl, Jacobi Operators and Completely Integrable Nonlinear Lattices.
Math. Surv. and Monographs, vol. 72, (2000) \\
6. Novikov ed. Soliton theory: the inverse problem method / "Nauka",
Moscow (1980). \\
7. Boutet de Monvel A. and Marchenko V. \,
The Cauchy problem for nonlinear Shr\"odinger equation with bounded
initial data. Mat. fiz., an., geom. (1997), v.4, ü~1/2, p.~3--45.
\end{flushleft}

\end{document}